\documentclass[11pt, reqno]{amsart}[11pt]
\usepackage{geometry}                
\usepackage{graphicx}
\usepackage{epstopdf}
\usepackage{amsmath,amssymb,amscd,mathrsfs,tikz-cd,mathtools,scalerel,latexsym}
\usepackage[perpage]{footmisc}
\usepackage{xcolor,todonotes}
\usepackage{pdfsync}
\usepackage[all]{xy}
\usepackage{tabularx}
\usepackage{enumitem}
\usepackage{autobreak}
\usepackage{breqn}
\usepackage{bm}
\usepackage{changes}[final]

\allowdisplaybreaks
\numberwithin{equation}{section}

\newtheorem{theorem}{Theorem}[section]
\newtheorem{remark}{Remark}[section]

\newcommand{\moda}{mod}

\newcommand{\fW}[1]{W\left(#1\right)}

\title{The second moment of Ramanujan sums}
\author{Hong Ziwei}
\address{School of Mathematics, Renmin University of China, Beijing, P.R. China}
\email{hongziwei@live.com}

\author{Zheng Zhiyong}
\address{School of Mathematics, 
Renmin University of China, Beijing, P.R. China}
\email{zhengzy@ruc.edu.cn}

\date{April 2023}

\begin{document}

\begin{abstract}
We study $C(x, y)$, the smoothed second moment of Ramanujan sums. Assuming the Riemann Hypothesis (RH), we establish an asymptotic formula for $C(x, y)$ with improved error term. Our analysis applies uniformly to the case where $x$ and $y$ are arbitrarily close. In particular, we provide a meaningful comparison with the work of \cite{TH} in case $y=2x^2$. The method relies on the use of smooth cutoff functions, which provide greater flexibility in contour shifting.
\end{abstract}

\maketitle



\section{Introduction}
The Ramanujan sum is a classical object in analytic number theory, defined by 
\begin{align}
    c_q(n)=\sum_{\substack{a\ \moda\ q\\ (a, q)=1}}e(\frac{na}{q}),
\end{align}
where $e(x)=e^{2\pi ix}$.

In \cite{TH}, T.H. Chan and A.V. Kumchev studied the second moment of the average of Ramanujan sums, defined by
\begin{equation}
    C_2(x, y)=\sum_{n\le y}\left(\sum_{q\le x}c_q(n)\right)^2.
\end{equation}

Using elementary methods, they obtained the asymptotic formula (see (1.2) in \cite{TH})
\begin{align}\label{1.2}
    C_2(x,y)=\frac{yx^2}{2\zeta(2)}+O(x^4+xy\log x).
\end{align}


Since their main interest lies in the asymptotic behavior of  $C_2(x,y)$ in the range $x^{1+\varepsilon}<y<x^{2+\varepsilon}$, Chan and Kumchev established the following theorem 
\begin{theorem}[T.H. Chan and A.V. Kumchev] 
    Let $x$ be a large real number, $y\ge x$, and $B>10$ be fixed.
    \begin{enumerate}
        \item If $y\ge x^2(\log x)^B$, then
        $$C_2(x, y)=\frac{yx^2}{2\zeta(2)}+O(yx^2(\log x)^{10}(x^{-1/2}+x^2y^{-1})).$$
        \item If $x(\log x)^{2B}\le y\le x^2(\log x)^B$, then
        $$C_2(x, y)=\frac{yx^2}{2\zeta(2)}(1+2\kappa(u))+O(yx^2(\log x)^{10}(x^{-1/2}+(y/x)^{-1/2})),$$
        where $u=\log(y/x^2)$, and $\kappa(u)$ is defined in \cite{TH}.
    \end{enumerate}
\end{theorem}
\begin{remark}
    This is the original version of the result of Chan and Kumchev. In the most recent version of \cite{TH}, the case $y\ge x^2(\log x)^B$ is replaced by the estimate in (\ref{1.2}).
\end{remark}
The above result improves the error term when $y=x^{2+\varepsilon}$. It is therefore natural to ask what can be said when $x$ and $y$ are of comparable size. \added{To address this question, we introduce a smooth version of the second moment considered by Chan and Kumchev}.

\subsection{Main result}

We begin by expressing $C_2(x, y)$ in a smooth form. Let $W:\mathbb{R}_{>0}\to\mathbb{R}$ be a smooth function such that $W(t)=1$ for $t\in(0,1)$, and $W(t)$ decays rapidly for $t\notin(0,1)$ (i.e. as $t\to 0^{+}$ and as $t\to +\infty$).
Moreover, for every integer $j\ge 0$ we have $W^{(j)}(t)\ll_{j} U^{\,j}$ for $t$ outside $(0,1)$, where $U\ge 2$ is a parameter. We then define the smooth second moment
\begin{equation}
    C(x, y)=C_2(x, y; W)=\sum_n\left(\sum_qc_q(n)\fW{\frac{q}{x}}\right)^2\fW{\frac{n}{y}}.
\end{equation}

The use of a smooth cutoff function allows us to exploit Mellin inversion and contour shifting more flexibly than in the nonsmooth setting.

Our main result is stated in the following theorem.
 \begin{theorem}
Assume the RH. Let $x$, $y$ be large real numbers with $y>x$, fix $W$. Let $\varepsilon>0$, $c_2\in (1, 2)$.
\begin{enumerate}
    \item If $x<y\le x^2$, then
    \begin{align}\label{y<x^2}
    C(x, y)=\frac{yx^2}{\zeta(2)}\widehat{W}(1)\frac{1}{2\pi i}\int\limits_{(c_2)}\widehat{W}(2-s_2)\widehat{W}(s_2)ds_2+O(y^{\varepsilon}x^{3-\varepsilon})+O(y^{\frac{3}{2}-\varepsilon}x^{1+2\varepsilon});
    \end{align}
    \item If $y> x^2$, then
    \begin{align}\label{y>x^2}
        C(x, y)=\frac{yx^2}{\zeta(2)}\widehat{W}(1)\frac{1}{2\pi i}\int\limits_{(c_2)}\widehat{W}(2-s_2)\widehat{W}(s_2)ds_2+O(yx^{\frac{1}{2}+\varepsilon})+O(y^{\varepsilon}x^{4-2\varepsilon}).
    \end{align}
\end{enumerate}
Here $\widehat{W}(s)$ is the Mellin transform of $W(t)$.
\end{theorem}
The Mellin transform of $W$ is defined by $\widehat{W}(s)=\int_{0}^{+\infty}W(t)t^{s-1}dt$. $\widehat{W}$ is analytic in the right half-plane $\sigma>0$, which allows coefficient $\widehat{W}(1)$ to be a computable constant. The integral in the equations has a magnitude of $\frac{1}{2}$. 

\added{
In fact, for $y>x$, we obtain a more precise asymptotic expansion of the form
\begin{theorem}
Assume the RH. Let $\varepsilon>0$, $c_1,c_2\in (1,2)$, large real number $y>x$, and $W$ be fixed. We have
\begin{align*}
    C(x,y)=\frac{yx^2}{\zeta(2)}\widehat{W}(1)\frac{1}{2\pi i}\int\limits_{(c_2)}\widehat{W}(2-s_2)\widehat{W}(s_2)ds_2+\frac{y^{2-c_1}x^{2c_1}}{\zeta(2)}\tilde{\kappa}(x,y,c_1)\\
    +O(yx^{\frac{1}{2}+\varepsilon})+O(y^{\frac{1}{2}+\varepsilon}x^{\frac{3}{2}+\varepsilon})+O(y^{\varepsilon}x^{3+\varepsilon})+O(y^{\frac{3}{4}}x^{1+\varepsilon}),
\end{align*}
where $\tilde{\kappa}(x,y,c_1)$ is a smooth analogue of $\kappa(u)$ in \cite{TH}, to be defined in Remark~\ref{K}. 
\end{theorem}}\added{
This expansion contains more detailed information about on the structure of $C(x,y)$. The finer decomposition plays a crucial role in the comparison of the results when $y$ is close to $x^2$.}
\added{
Since our primary interest lies in the range $x^{1+\varepsilon}<y<x^{2+\varepsilon}$, we state our main result in the simplified form of Theorem~1.2. A detailed analysis of the secondary main term (the choice of parameter $c_1$) is carried out in Section~3.} 

Finally, we note that there has been growing interest in the study of mean values of Ramanujan sums over number fields, leading to a number of recent developments (see, for instance, \cite{ref:Nowak, Nowak, ref:JCG, ref:SC}). The method developed in this paper can be extended to the setting of number fields and is expected to yield corresponding improvements in that context. For the sake of clarity, and following the approach of Chan and Kumchev, we restrict our discussion here to the ring of integers, which already illustrates the ideas of the argument.

\subsection{Comparison}
\added{There are three main differences between the result of Chan and Kumchev and the present work.}

\added{
First, our method yields an asymptotic formula uniformly in ranges where $y$ may be as close to $x$ as $y/x\to 1$.
This relies on the smooth cutoff (see Remark~\ref{shifting}), which enables a larger contour shift than is available with a sharp cutoff.} 

\added{
Second, in certain ranges where $y$ is large, our estimate improves upon those obtained by Chan and Kumchev. For instance, when $y=2x^2$, both results yield a main term of order $x^4$. Combining the estimate \eqref{eq:C_3} for $C_3$ with the bounds for other terms, 
we obtain an error term of order $x^{3+\varepsilon}$, whereas the result of Chan and Kumchev yields an error term of size $x^{\frac{7}{2}}(\log x)^{10}$.} 


\added{
Furthermore, Chan?s method relies on the computation of integrals over a zero-free region of $\zeta(s)$. Under RH, all zeros of $\zeta(s)$ are on the line $Re(s)=\frac{1}{2}$ which allows larger zero-free region and sharper bounds of $C(x,y)$. 
}

 



\section{Preliminary}
 \subsection{Mellin inversion}
Let $W$ be a smooth, rapidly decaying function. By Mellin inversion, we may write
 \begin{equation}\label{Mellin}
     \fW{t}=\frac{1}{2\pi i}\int\limits_{(c)}\widehat{W}(s)t^{-s}ds,
 \end{equation}
 where $\widehat{W}(s)$ denotes the Mellin transform  of $W$, defined by 
 \begin{equation}\label{Mi}
 \widehat{W}(s)=\int_{0}^{+\infty}W(t)t^{s-1}dt.
 \end{equation}
Here, $c>1$ is fixed, and $\widehat{W}(s)$ is analytic for $Re(s)>0$. Moreover, for any positive integer $E$, we have the bound 
 \begin{equation}\label{estimateW}
     \widehat{W}(s)\ll \frac{1}{|s|(1+|s|)^E}.
 \end{equation}
 For further details on $W$ and $\widehat{W}$, we refer to \cite{Gao}. In the present work, the precise form of $W$ is immaterial; only the rapid decay of $\widehat{W}$ will be used.
 
 \subsection{Ramanujan's identity}
Applying Mellin inversion \eqref{Mellin} together with Ramanujan's identity (see (3.1) in \cite{TH}), we may express $C(x,y)$ as
 \begin{align*}
     C(x, y)&=\sum_n\left(\frac{1}{2\pi i}\int\limits_{(c)}\sum_q\frac{c_q(n)}{q^s}x^s\widehat{W}(s)ds\right)^2\fW{\frac{n}{y}}\notag\\
     &=\sum_n\left(\frac{1}{2\pi i}\int\limits_{(c)}\frac{\sigma_{1-s}(n)}{\zeta(s)}x^s\widehat{W}(s)ds\right)^2\fW{\frac{n}{y}}\notag\\
     &=\sum_n\left(\frac{1}{2\pi i}\right)^2\int\limits_{(c_1)}\int\limits_{(c_2)}\frac{\sigma_{1-s_1}(n)\sigma_{1-s_2}(n)}{\zeta(s_1)\zeta(s_2)}x^{s_1+s_2}\widehat{W}(s_1)\widehat{W}(s_2)ds_1ds_2\fW{\frac{n}{y}}
 \end{align*}
 where $2>c_2>c_1>1$ are real numbers close to $1$. Combining this with equation (3.6) of \cite{TH}, we may rewrite the above expression as
 \begin{align}
     C(x, y)&=\left(\frac{1}{2\pi i}\right)^2\int\limits_{(c_1)}\int\limits_{(c_2)}\frac{\sum\limits_n\sigma_{1-s_1}(n)\sigma_{1-s_2}(n)\fW{\frac{n}{y}}}{\zeta(s_1)\zeta(s_2)}x^{s_1+s_2}\widehat{W}(s_1)\widehat{W}(s_2)ds_1ds_2\notag\\
     &=\left(\frac{1}{2\pi i}\right)^2\int\limits_{(c_1)}\int\limits_{(c_2)}H(s_1, s_2; y)\frac{x^{s_1+s_2}}{\zeta(s_1)\zeta(s_2)}\widehat{W}(s_1)\widehat{W}(s_2)ds_1ds_2,
 \end{align}
 where
 \begin{align}\label{H}
     H(s_1, s_2; y)=\frac{1}{2\pi i}\int\limits_{(c_3)}\frac{\zeta(s_3)\zeta(s_3+s_1-1)\zeta(s_3+s_2-1)\zeta(s_3+s_1+s_2-2)}{\zeta(2s_3+s_1+s_2-2)}\widehat{W}(s_3)y^{s_3}ds_3.
 \end{align}
 Here $c_3$ is chosen so that $c_3>\max\{1, 2-Re(s_1), 2-Re(s_2), 3-Re(s_1+s_2)\}$, which ensures absolute convergence of the integral. In what follows, it will be sufficient to take $c_3=3/2$. Further details on the convergence regions of the Riemann zeta function $\zeta(s)$ may be found in \cite{ref:Multi}.
 
 \section{The second moment}
 \subsection{Expression for $C(x,y)$}
  The integrand in \eqref{H} has poles at $s_3=1$, $s_3=2-s_1$, $s_3=2-s_2$ and $s_3=3-s_1-s_2$ with corresponding residues $$\frac{\zeta(s_1)\zeta(s_2)\zeta(s_1+s_2-1)}{\zeta(s_1+s_2)}y\widehat{W}(1),$$
 $$\frac{\zeta(2-s_2)\zeta(1-s_2+s_1)\zeta(s_1)}{\zeta(2+s_1-s_2)}y^{2-s_2}\widehat{W}(2-s_2),$$
 $$\frac{\zeta(2-s_1)\zeta(1-s_1+s_2)\zeta(s_2)}{\zeta(2+s_2-s_1)}y^{2-s_1}\widehat{W}(2-s_1)$$
 and 
 $$\frac{\zeta(3-s_1-s_2)\zeta(2-s_1)\zeta(2-s_2)}{\zeta(4-(s_1+s_2))}y^{3-s_1-s_2}\widehat{W}(3-s_1-s_2)$$
 respectively.
 Assuming the Riemann Hypothesis, we shift the contour of integration from $(c_3)$ to $(c_3')$, with $c_3'<1$, and collecting the residues at the poles listed above, we may write 
 \begin{align*}
     H(s_1, s_2; y)&=\frac{\zeta(s_1)\zeta(s_2)\zeta(s_1+s_2-1)}{\zeta(s_1+s_2)}y\widehat{W}(1)\\
     &+\frac{\zeta(2-s_2)\zeta(1-s_2+s_1)\zeta(s_1)}{\zeta(2+s_1-s_2)}y^{2-s_2}\widehat{W}(2-s_2)\\
     &+\frac{\zeta(2-s_1)\zeta(1-s_1+s_2)\zeta(s_2)}{\zeta(2+s_2-s_1)}y^{2-s_1}\widehat{W}(2-s_1)\\
     &+\frac{\zeta(3-s_1-s_2)\zeta(2-s_1)\zeta(2-s_2)}{\zeta(4-(s_1+s_2))}y^{3-s_1-s_2}\widehat{W}(3-s_1-s_2)\\
     &+\frac{1}{2\pi i}\int\limits_{(c_3')}\frac{\zeta(s_3)\zeta(s_3+s_1-1)\zeta(s_3+s_2-1)\zeta(s_3+s_1+s_2-2)}{\zeta(2s_3+s_1+s_2-2)}\widehat{W}(s_3)y^{s_3}ds_3.
 \end{align*}

 Substituting this decomposition into $C(x, y)$, we obtain
 \begin{align}
     C(x, y)=C_1+C_2+C_3+C_4+R,
 \end{align}
 where
 \begin{align*} 
     C_1=\left(\frac{1}{2\pi i}\right)^2\int\limits_{(c_1)}\int\limits_{(c_2)}\frac{\zeta(s_1+s_2-1)}{\zeta(s_1+s_2)}y\widehat{W}(1)x^{s_1+s_2}\widehat{W}(s_1)\widehat{W}(s_2)ds_1ds_2,
 \end{align*}
 \begin{align*} 
     C_2=\left(\frac{1}{2\pi i}\right)^2\int\limits_{(c_1)}\int\limits_{(c_2)}\frac{\zeta(2-s_2)\zeta(1-s_2+s_1)}{\zeta(2+s_1-s_2)\zeta(s_2)}y^{2-s_2}\widehat{W}(2-s_2)x^{s_1+s_2}\widehat{W}(s_1)\widehat{W}(s_2)ds_1ds_2,
 \end{align*}
 \begin{align*}  
     C_3=\left(\frac{1}{2\pi i}\right)^2\int\limits_{(c_1)}\int\limits_{(c_2)}\frac{\zeta(2-s_1)\zeta(1-s_1+s_2)}{\zeta(2+s_2-s_1)\zeta(s_1)}y^{2-s_1}\widehat{W}(2-s_1)x^{s_1+s_2}\widehat{W}(s_1)\widehat{W}(s_2)ds_1ds_2,
 \end{align*}
 \begin{align*} 
     C_4=\left(\frac{1}{2\pi i}\right)^2\int\limits_{(c_1)}\int\limits_{(c_2)}\frac{\zeta(3-s_1-s_2)\zeta(2-s_1)\zeta(2-s_2)}{\zeta(4-(s_1+s_2))\zeta(s_1)\zeta(s_2)}y^{3-s_2-s_1}x^{s_1+s_2}\\
     \times\widehat{W}(3-s_1-s_2)\widehat{W}(s_1)\widehat{W}(s_2)ds_1ds_2
 \end{align*}
 and
 \begin{align*} 
     R=\left(\frac{1}{2\pi i}\right)^3\int\limits_{(c_1)}\int\limits_{(c_2)}\int\limits_{(c_3')}\frac{\zeta(s_3)\zeta(s_3+s_1-1)\zeta(s_3+s_2-1)\zeta(s_3+s_1+s_2-2)}{\zeta(2s_3+s_1+s_2-2)\zeta(s_1)\zeta(s_2)}y^{s_3}x^{s_1+s_2}\\
     \times\widehat{W}(s_3)\widehat{W}(s_1)\widehat{W}(s_2)ds_1ds_2ds_3. 
     \end{align*}

 \subsection{The first main term in $C(x, y)$}
 We begin by analyzing the contribution of $C_1$. The function $\zeta(s_1+s_2-1)$ has a simple pole at $s_2=2-s_1$ with residue 1. Assuming the Riemann Hypothesis, we shift the contour from $(c_2)$ to $(c_2')$ with $c_1+c_2'-1<1$. Taking the residue at $s_2=2-s_1$, we obtain the main term. More precisely, we may write 
 \begin{align*}
     C_1=&\frac{yx^2}{\zeta(2)}\widehat{W}(1)\frac{1}{2\pi i}\int\limits_{(c_1)}\widehat{W}(2-s_1)\widehat{W}(s_1)ds_1\\
     &+\left(\frac{1}{2\pi i}\right)^2y\widehat{W}(1)\int\limits_{(c_1)}\int\limits_{(c_2')}\frac{\zeta(s_1+s_2-1)}{\zeta(s_1+s_2)}x^{s_1+s_2}\widehat{W}(s_1)\widehat{W}(s_2)ds_1ds_2.
 \end{align*}
 The size of the second term in the above equation is governed by the factor $x^{c_1 + c_2'}$. We choose $c_1$ and $c_2'$ so as to minimize $c_1+c_2'$ subject to the constraints $c_1+c_2'>\frac{1}{2}$ and $c_1, c_2'>0$. Taking $c_1=c_2'=\frac{1}{4}+\frac{\varepsilon}{2}$, we obtain
 \begin{align*}
     &\ll yx^{\frac{1}{2}+\varepsilon}\int^{+\infty}_{-\infty}\int^{+\infty}_{-\infty}\frac{|\zeta(-\frac{1}{2}+\varepsilon+(t_1+t_2)i)|}{|\zeta(\frac{1}{2}+\varepsilon+(t_1+t_2)i)|}\widehat{W}(\frac{1}{4}+\frac{\varepsilon}{2}+t_1i)\widehat{W}(\frac{1}{4}+\frac{\varepsilon}{2}+t_2i)dt_1dt_2.
 \end{align*}
 For $|\frac{1}{\zeta(s)}|$, we use the bound $(14.2.6)$ of \cite{Titchmarsh}. By the functional equation and Stirling's formula, one obtains the standard bound 
 \begin{align}
        |\zeta(s)|\ll\left(\frac{|t|}{2\pi}\right)^{\frac{1}{2}-\sigma}\left(1+\frac{1}{|t|}\right)|\zeta(1-\sigma)|
 \end{align} 
 valid for $-1\le \sigma\le\frac{1}{2}$. Together with the rapid decay of $\widehat{W}$ this shows that the above integral is convergent. Consequently, we arrive at
 \begin{align}\label{C1}
     C_1=\frac{yx^2}{\zeta(2)}\widehat{W}(1)\frac{1}{2\pi i}\int\limits_{(c_2)}\widehat{W}(2-s_2)\widehat{W}(s_2)ds_2+O(yx^{\frac{1}{2}+\varepsilon}).
 \end{align}
 
 It's worth noting that \eqref{Mi} implies $\widehat{W}(1)=1+O(\frac{1}{U})$ and $\frac{1}{2\pi i}\int\limits_{(c_2)}\widehat{W}(2-s_2)\widehat{W}(s_2)ds_2=\frac{1}{2}+O(\frac{1}{U})$. 

\subsection{Discussion of the size of $\frac{x^2}{y}$}\label{discussion}
 We now focus on the analysis of the term $C_3$. Throughout this section, we assume $y>x$ and that $2>c_2>c_1>1$. By definition, 
 \begin{equation*}
     C_3=\left(\frac{1}{2\pi i}\right)^2\int\limits_{(c_1)}\int\limits_{(c_2)}\frac{\zeta(2-s_1)\zeta(1-s_1+s_2)}{\zeta(2+s_2-s_1)\zeta(s_1)}y^{2}\widehat{W}(2-s_1)x^{s_2}\left(\frac{x}{y}\right)^{s_1}\widehat{W}(s_1)\widehat{W}(s_2)ds_1ds_2.
 \end{equation*}
 The function $\zeta(1-s_1+s_2)$ has a simple pole at $s_2=s_1$ with residue $1$. Shifting the contour from $(c_2)$ to $(c_2')$, where $1-c_1+c_2'<1$ and $c_2'<c_1$, we obtain
 \begin{align*}
     C_3=&\frac{y^2}{\zeta(2)}\frac{1}{2\pi i}\int\limits_{(c_1)}\frac{\zeta(2-s_1)}{\zeta(s_1)}\left(\frac{x^2}{y}\right)^{s_1}\widehat{W}(2-s_1)\widehat{W}(s_1)^2ds_1\\
     &+\frac{1}{(2\pi i)^3}\int\limits_{(c_1)}\int\limits_{(c_2')}\frac{\zeta(2-s_1)\zeta(1-s_1+s_2)}{\zeta(2+s_2-s_1)\zeta(s_1)}y^{2}\widehat{W}(2-s_1)x^{s_2}\left(\frac{x}{y}\right)^{s_1}\widehat{W}(s_1)\widehat{W}(s_2)ds_1ds_2.
 \end{align*}

The size of the second term is governed by $Re(s_1)$. Shifting the contour from $(c_1)$ to $(\frac{3}{2}-\frac{\varepsilon}{2})$ and $(c_2')$ to $(\varepsilon)$, and arguing as in the estimate of $C_1$, we obtain
 \begin{align}\label{eq:C_3}
     C_3=\frac{y^2}{\zeta(2)}\frac{1}{2\pi i}\int\limits_{(c_1)}\frac{\zeta(2-s_1)}{\zeta(s_1)}\left(\frac{x^2}{y}\right)^{s_1}\widehat{W}(2-s_1)\widehat{W}(s_1)^2ds_1+O(y^{\frac{1}{2}+\varepsilon}x^{\frac{3}{2}+\varepsilon}).
 \end{align} \added{
 \begin{remark}\label{K}
     Define
     \begin{align*}
         \tilde{\kappa}(x,y,c_1)=\frac{1}{2\pi i}\int_{(c_1)}\frac{\zeta(2-s_1)}{\zeta(s_1)}\left(\frac{x^2}{y}\right)^{Im(s_1)}\widehat{W}(2-s_1)\widehat{W}(s_1)^2ds_1.
     \end{align*}
     Then the contribution $C_3$ can be written as
     \begin{align*}
         C_3=\frac{y^{2-c_1}x^{2c_1}}{\zeta(2)}\tilde{\kappa}(x,y,c_1)+O(y^{\frac{1}{2}+\varepsilon}x^{\frac{3}{2}+\varepsilon}).
     \end{align*}
     The function $\tilde{\kappa}(x,y,c_1)$ may be viewed as a smooth version of the function $\kappa(u)$ appearing in the work of Chan and Kumchev. 
 \end{remark}}
 The first term in $C_3$ is influenced by the magnitude of $\frac{x^2}{y}$.

 \textit{Case 1.} $y=x^2$
 
 In this case, the first term of $C_3$ simplifies to
 \begin{align*}
     \frac{y^2}{\zeta(2)}\frac{1}{2\pi i}\int\limits_{(c_1)}\frac{\zeta(2-s_1)}{\zeta(s_1)}\widehat{W}(2-s_1)\widehat{W}(s_1)^2ds_1,
 \end{align*}
 where $c_1\in (\frac{1}{2}, 2)$. Since the integrand is analytic in this strip, the integral is a constant. Consequently,
 \begin{align*}
     C_3\ll y^2+y^{\frac{1}{2}+\varepsilon}x^{\frac{3}{2}+\varepsilon}\ll y^2.
 \end{align*}

 \textit{Case 2.} $x<y<x^2$
 
 In this case, we have $\frac{x^2}{y}>1$. Shifting the contour from $(c_1)$ to $(\frac{1}{2}+\varepsilon)$, and assuming the Riemann Hypothesis, we obtain
 \begin{align*}
     C_3=&\frac{y^2}{\zeta(2)}\frac{1}{2\pi i}\int\limits_{(\frac{1}{2}+\varepsilon)}\frac{\zeta(2-s_1)}{\zeta(s_1)}\left(\frac{x^2}{y}\right)^{s_1}\widehat{W}(2-s_1)\widehat{W}(s_1)^2ds_1+O(y^{\frac{1}{2}+\varepsilon}x^{\frac{3}{2}+\varepsilon})\\
     \ll &y^{\frac{3}{2}-\varepsilon}x^{1+\varepsilon}+y^{\frac{1}{2}+\varepsilon}x^{\frac{3}{2}+\varepsilon}\\
     \ll &y^{\frac{3}{2}-\varepsilon}x^{1+\varepsilon}.
 \end{align*}

 \textit{Case 3.} $y>x^2$
 
 In this case, $\frac{x^2}{y}<1$. Shifting the contour from $(c_1)$ to $(2-\varepsilon)$, we deduce
 \begin{align*}
     C_3\ll &y^{\varepsilon}x^{4-2\varepsilon}+y^{\frac{1}{2}+\varepsilon}x^{\frac{3}{2}+\varepsilon}.
 \end{align*}
  \begin{remark}
      The above discussion implicitly involves a comparison between $y$ and $x^3$. Since our primary interest lies on the range $x<y<x^{2+\varepsilon}$, we do not pursue further refinements in the remaining ranges.
  \end{remark} 

 Collecting the above estimates, we obtain
 \begin{equation}\label{C3}
     C_3\ll\left\{\begin{array}{ll}
         y^{\frac{3}{2}-\varepsilon}x^{1+\varepsilon}+y^{\frac{1}{2}+\varepsilon}x^{\frac{3}{2}+\varepsilon}, & x<y\le x^2;\\
         y^{\varepsilon}x^{4-2\varepsilon}+y^{\frac{1}{2}+\varepsilon}x^{\frac{3}{2}+\varepsilon}, & y>x^2.
     \end{array}
     \right.
 \end{equation}

 \subsection{Error terms}
 We now estimate the remaining terms $C_2$, $C_4$, and $R$. We begin with the term $C_2$, which is given by
 $$C_2=\left(\frac{1}{2\pi i}\right)^2\int\limits_{(c_1)}\int\limits_{(c_2)}\frac{\zeta(2-s_2)\zeta(1-s_2+s_1)}{\zeta(2+s_1-s_2)\zeta(s_2)}y^2\widehat{W}(2-s_2)x^{s_1}\left(\frac{x}{y}\right)^{s_2}\widehat{W}(s_1)\widehat{W}(s_2)ds_1ds_2.$$
 Our aim is to choose the contours so that $Re(s_1)$ is small and $Re(s_2)$ is large. Assuming the Riemann Hypothesis, the factors $\zeta(2+s_1-s_2)$ and $\zeta(s_2)$ impose the constraints $Re(2+s_1-s_2)>\frac{1}{2}$ and $Re(s_2)>\frac{1}{2}$. Shifting the contour from $(c_2)$ to $(\frac{3}{2}-\varepsilon/2)$ and from $(c_1)$ to $(\varepsilon)$, we obtain 
 \begin{align}\label{C2}
     C_2\ll y^{\frac{1}{2}+\varepsilon}x^{\frac{3}{2}+\varepsilon}.
 \end{align}
 Arguing similarly and taking $c_1=c_2=\frac{3}{2}-\frac{\varepsilon}{2}$, we deduce that
 \begin{align}\label{C4}
     C_4\ll y^{\varepsilon}x^{3-\varepsilon}.
 \end{align}
\added{
\begin{remark}\label{shifting}
 When $y=x^{1+\varepsilon}$, $C_4$ gives the upper bound of the error term. Chan and Kumchev employ an unsmooth cutoff function, which leads to a restriction of the form $Re(3-(s_1+s_2))=\frac{1}{2}$ in the contour shifting (see (4.15) in \cite{TH}). By contrast, the use of a smooth cutoff function in our approach allows us to relax the constraint to $Re(3-(s_1+s_2))>0$ (see the expression of $C_4$ in Section~3.1). The additional flexibility is crucial when $y$ is very close to $x$.
 \end{remark}}
 
 Finally, we estimate the remainder term $R$. The integrand imposes the constraints $Re(s_1)>\frac{1}{2}$, $Re(s_2)>\frac{1}{2}$ and $Re(2s_3+s_1+s_2)>\frac{5}{2}$. To minimize $Re(s_3)$ and $Re(s_1+s_2)$, we choose $c_3'=\frac{3}{4}$ and $c_1=c_2=\frac{1}{2}+\frac{\varepsilon}{2}$. With this choice, we obtain 
 \begin{align}\label{R}
     R\ll y^{\frac{3}{4}}x^{1+\varepsilon}.
 \end{align}

 \section{Conclusion}

 Combining the estimates for the terms $C_1$, $C_2$, $C_3$, $C_4$, and $R$ as given in equations \eqref{C1}, \eqref{C2}, \eqref{C3}, \eqref{C4}, and \eqref{R}, respectively, we complete the proof of Theorem~1.2.
 
  \section*{Acknowledgement}
The authors thank Professor P. Gao for suggesting the subject of this study and for providing invaluable guidance and insightful recommendations throughout the research process.

\bibliography{biblio}
\bibliographystyle{plain}

\end{document}